\documentclass[intlimits,12pt]{amsart}
\usepackage{hyperref}
\usepackage[all]{xy}
\usepackage{amssymb}

\newcommand{\vp}{\varphi}
\newcommand{\begeq}{\begin{equation}}

\newcommand{\dra}{\rangle\rangle}
\newcommand{\dla}{\langle\langle}

\newcommand{\zzendeq}{\end{equation}}
\newcommand{\begmat}{\left(\begin{array}{cc}}
\newcommand{\zzendmat}{\end{array}\right)}
\newcommand{\R}{\bold{R}}
\newcommand{\Unn}{U_{\nu_1,\nu_2}}
\newcommand{\Q}{\bold{Q}}

\newcommand{\0}{\boldsymbol{0}}
\newcommand{\Z}{\bold{Z}}
\newcommand{\N}{\bold{N}}
\newcommand{\C}{\bold{C}}

\newcommand{\bQ}{\mathbf{Q}}
\newcommand{\bZ}{\mathbf{Z}}

\newcommand{\cA}{{\mathcal A}}

\newcommand{\A}{\bold{A}}

\newcommand{\Sum}{\sum\limits}

\newcommand{\tr}{\text{tr}}

\newcommand{\Endo}{\text{End}}
\newtheorem{theorem}{Theorem}[section]
\newtheorem{lemma}[theorem]{Lemma}

\newtheorem{corollary}[theorem]{Corollary}
\newtheorem{proposition}[theorem]{Proposition}
\newtheorem{remark}[theorem]{Remark}
\numberwithin{equation}{section}

\newcommand{\tdiff}{\tilde{{\mathcal D}}^\gamma_{2,\alpha_1,\alpha_2}}
\begin{document}
\title[Global Gross-Prasad conjecture for Yoshida liftings]{On the global Gross-Prasad conjecture for Yoshida liftings}
\author[S. B\"ocherer \and M. Furusawa \and  R. Schulze-Pillot]{Siegfried B\"ocherer \and Masaaki Furusawa \and  Rainer Schulze-Pillot} 
\thanks{Part of this work was done during a stay of all three authors
  at the IAS, Princeton, supported by von Neumann Fund, Bell
  Fund and Bankers
  Trust Fund. Furusawa was also supported by Sumitomo
  Foundation. We thank the IAS for its hospitality.} 
\thanks {S. B\"ocherer and R. Schulze-Pillot thank D. Prasad and the Harish Chandra
  Research Institute, Allahabad, India for their
  hospitality. Schulze-Pillot's visit to HCRI was also supported by DFG}  
\dedicatory{Prof. J. Shalika to his 60th birthday}
\begin{abstract}
We restrict a Siegel modular cusp form of degree $2$ and square free
level that is a Yoshida lifting (a lifting from the orthogonal group
of a definite quaternion algebra) 
to the embedded product of two half planes and compute the Petersson
product against the product of two elliptic cuspidal Hecke eigenforms.
The square of this integral can be explicitly expressed in terms of
the central critical value of an $L$-function attached to the
situation. The result is related to a conjecture of Gross and Prasad about
restrictions of automorphic representations of special orthogonal groups

\end{abstract}

\maketitle

\section*{Introduction}
In two articles in the Canadian Journal
\cite{gropra1,gropra2}, B. Gross and D. Prasad 
proclaimed a global
conjecture concerning the decomposition of an automorphic
representation of an adelic special orthogonal group $G_1$
upon restriction to an embedded orthogonal
group $G_2$ of a
quadratic space in smaller dimension
and also its local counterpart.
In the local
situation, one can summarize the conjecture by saying that
the occurrence of $\pi_2$ in the restriction of $\pi_1$
depends on the $\epsilon$-factor attached to the
representation
$\pi_1\otimes \pi_2;$
in the global situation, assuming the existence of
the local nontrivial invariant functional at all places
and its nonvanishing on the spherical vector at almost
all unramified places,
one considers a specific linear functional given by a period integral.
This period integral is then conjectured to give a nontrivial
functional if and only if
the central critical value
of the  $L$-function attached to $\pi_1\otimes\pi_2$
is nonzero.
In particular in the
case when $G_1$ is the group of an $n$-dimensional nondegenerate 
quadratic space $V$ and $G_2$ is
the group of an $(n-1)$-dimensional subspace $W$ of $V$,
they showed that in low dimensions ($n\le 4$) known results
can be interpreted as evidence for this conjecture, using
the well known isomorphisms for orthogonal groups in low
dimensions.  

The case $n=5$ has been treated in the local situation by Prasad
\cite{prasad_seesaw}; it can also be reinterpreted using
these isomorphisms: The split special orthogonal group in
dimension $5$ is isomorphic to the projective symplectic
similitude group $\mathrm{PGSp}_2,$ and the spin group of the
$4$-dimensional split orthogonal group is $\mathrm{SL}_2 \times
\mathrm{SL}_2.$ Prasad then showed that for forms on
$\mathrm{PGSp}_2$ that are lifts from the orthogonal group of a
4-dimensional space, the situation can be
understood in terms of the seesaw dual reductive pair (in Kudla's
sense)
$$\xymatrix{GSp_2 \ar@{-}[d]\ar@{-}[dr]&GO(4)\ar@{-}[d]\times GO(4)\\
  G(SL_2 \times SL_2) \ar@{-}[ur]&GO(4)}.$$

In classical terms, the analogous global question leads to the problem
to determine those 
pairs of cuspidal elliptic modular forms (eigenforms of almost all
Hecke operators) that can occur as summands if one decomposes the
restriction of a  cuspidal 
Siegel modular form $F$ of degree $2$ (that is an eigenform of almost
all  Hecke operators) to the diagonally embedded product of two upper 
half planes into a sum of (products of) eigenforms of almost all Hecke operators
(for a discussion of the problems that arise in translating the
representation theoretic statement into a classical statement see below
and Remark \ref{rep_remark}. 
One can also rephrase
this as the problem of calculating the period integral 
$$\int_{\left(\Gamma \backslash 
{\bf H}\right)\times 
\left(\Gamma \backslash {\bf H}\right)}
F\left(\begin{pmatrix}
    z_1&0\\0&z_2\end{pmatrix}\right)
\overline{f_1(z_1)}\overline{f_2(z_2)}d^*z_1d^*z_2$$ for two
elliptic Hecke eigenforms $f_1,f_2.$ 
The $L$-function
that should occur then according to the conjecture of Gross
and Prasad is the degree $16$ $L$-function associated to the
tensor product of the $4$-dimensional representation of the
$L$-group $\mathrm{Spin}(4)$ of $\mathrm{SO}(5)$ with the 
two $2$-dimensional
representation associated to two copies of $\mathrm{SO}(3)$ (due to
the decomposition of $\mathrm{SO}(4)$ or rather its covering group
mentioned above). We denote this $L$-function
as $L({\mathrm
  Spin}(F),f_1,f_2,s).$

In this reformulation it is natural to go beyond the original
question of nonvanishing and to try and get an explicit formula
connecting the $L$-value in question with the period integral.
In general, it seems rather difficult to calculate the
period integral as above, since little is known about
the restriction of Siegel modular forms to the diagonally
embedded product of two upper half planes.
One should also point out that
an integral representation for the
degree $16$ $L$-function
in question is not known yet.
However, for
theta series of quadratic forms the restriction to the
diagonal of a degree two theta series becomes simply the
product of the degree one theta series in the variables
$z_1,z_2,$ which allows one to get a calculation started.
Thus here we only
consider Siegel modular forms (of trivial
character) that arise as linear
combinations of theta series of quaternary quadratic
forms. Such Siegel modular forms, if they are eigenforms,
are called Yoshida liftings,
attached to a pair of elliptic cusp
forms or, equivalently, to a pair of automorphic forms on
the multiplicative group of an adelic quaternion algebra.
These liftings have been investigated in
\cite{yoshida1,bs-nagoya1,bs-nagoya2}, and 
the connection between
the trilinear forms on the spaces of
automorphic forms on the multiplicative group of an adelic
quaternion algebra and  the triple product $L$-function has
been investigated in \cite{HK,gross-kudla,bs-triple}. If one
combines these results and applies them to the present
situation,
it turns out that the period integral in question
can indeed be explicitly calculated in terms of the central
critical value of the $L$-function mentioned; in the case of
a Yoshida-lifting attached to the pair $h_1,h_2$ of
elliptic cusp forms,
this $L$-function is seen to split into the
product
$L(h_1,f_1,f_2,s)L(h_2,f_1,f_2,s),$
so that the central critical value becomes the product of
the central critical values of these two triple product
$L$-functions.
We prove a formula that expresses the square of the period
integral explicitly as the product of these two central
critical values, multiplied by an explicitly known non-zero
factor. We reformulate  the obtained
identity in a way 
which
makes sense as well for an arbitrary Siegel modular form
$F$ 
in terms of 
the original $L({\mathrm Spin}(F),f_1,f_2,s)$, in place of 
 $L(h_1,f_1,f_2,s)L(h_2,f_1,f_2,s)$,
hoping  that such an identity indeed holds
for any Siegel modular form $F$.
At present we cannot prove it 
except for the case when $F$ is a Siegel or Klingen
Eisenstein series of level $1.$ In the case when $F$ is the
Saito-Kurokawa lifting
of an elliptic  Hecke eigenform $h$, one
sees easily that the period 
integral is zero unless one has $f_1=f_2=f$; in this case the
period integral can be transformed into the
Petersson inner product of the restriction of
the first Fourier Jacobi coefficient of $F$ to the
upper half plane with $f$ and then leads us  to a
conjectural identity for the square of this Petersson inner
product with the central critical value of $L(h,f,f;s).$  

Our calculation leaves the question open whether it can
happen that the period integral vanishes for the classical
modular forms considered but is non-zero for other functions
in the same adelic representation space. Viewed locally,
this amounts to the question whether an invariant nontrivial
linear 
functional on the local representation space is necessarily 
non-zero at the given vector. At the infinite place we can
exhibit such a vector (depending on the weights given) by
applying a suitable differential operator to the 
Siegel modular form considered. At the finite places not
dividing the level, it comes down to the question whether
(for an unramified representation) a 
nontrivial invariant linear functional is necessarily
non-zero at the spherical (or class 1) vector invariant under
the maximal compact subgroup. This is generally expected, at
least for generic representations. We intend to come back to
this question in future work.

We also investigate the situation where the pair
$f_1,f_2$ and the product ${\mathbf H}\times {\mathbf H}$ are replaced by a
Hilbert modular form and the modular embedding of a Hilbert
modular surface; in terms of the Gross-Prasad conjecture this
amounts to replacing the split orthogonal group of a
$4$-dimensional space from above
by a non split (but quasisplit) orthogonal group that is
split at infinity. It turns out that one gets an analogous
result; we prove this only in the simplest case when
all
modular forms involved have weight $2,$ the class number of
the quadratic field involved is $1$ and the order in a
quaternion algebra belonging to the situation is a maximal
order.  The proof for the general case should be
possible in an analogous
manner. 

\section{Yoshida liftings and their restriction to the diagonal}
For generalities on Siegel modular forms we refer to [Fre1].
For a symplectic matrix $M=
\left(\begin{array}{cc}
A & B \\
C & D       \end{array} \right)\in 
GSp_n({\R})$ 
(with $n\times
n$-blocks $A,B,C,D$)
we denote by $(M,Z)\mapsto M<Z>=( A Z+ B)( C Z+ D)^{-1}$  the usual action
of the group $G^{+}Sp(n,\R)$  of proper symplectic similitudes
on Siegel's upper half space ${\mathbf H}_{n}$.    \\
We shall mainly be concerned with Siegel modular forms for
congruence subgroups of type
$$\Gamma_{0}^{(n)}(N) = \left\{ \left(\begin{array}{cc}
 A &  B \\
 C &  D       \end{array} \right) \in Sp(n,\bZ) \mid  C \equiv 0 \bmod N
\right\}         .$$
The space of Siegel modular forms (and cusp forms respectively)
of degree $n$ and weight $k$ for $\Gamma_{0}^{(n)}(N)$
will be denoted by $M_{n}^{k}(N)$ $(S_{n}^{k}(N)),$ for a vector
valued modular form transforming according to the representation
$\rho$ the weight $k$ above should be replaced by $\rho.$ By $<~,~>$
we denote the Petersson scalar product.

We recall from \cite{bs-nagoya1,bs-nagoya2,bs-nachrichten} some
notations concerning the Yoshida-liftings whose restrictions we are
going to study in this article. For details we refer to the cited articles.
We consider a definite quaternion algebra $D$ over $\Q$ and 
an Eichler order $R$ 
of square free level $N$ in it
and decompose $N$  as
$N=N_1N_2$ where $N_1$ is the product of the primes that are
ramified in $D$.
On $D$ we have the involution $x\mapsto\overline x$, the
(reduced) trace
tr$(x)=x+\overline x$ and the (reduced) norm $n(x)=x\overline x$.

The group
of proper similitudes of the quadratic form $q(x)=n(x)$  
on $D$ is isomorphic to
$(D^\times\times D^\times)/Z(D^\times)$ (as algebraic group)
via 
$$(x_1,x_2)\mapsto \sigma_{x_1,x_2}\text{ with }
\sigma_{x_1,x_2}(y)  =  
x_1 y x_2^{-1} ,$$
the special orthogonal group is then the image of
$$\{(x_1,x_2)\in D^\times\times D^\times\mid n(x_1)= n(x_2)\} .$$
We denote by $H$ the orthogonal group of $(D,n)$ and by $H^+$ 
the special orthogonal
group.

For $\nu \in \N$ let $U_{\nu}^{(0)} $ be the space of homogeneous
harmonic polynomials of degree $\nu$ on $\R^3$ and view $P \in 
U_{\nu}^{(0)}$ as a polynomial on 
$$D_\infty^{(0)}= \{ x \in D_\infty
\vert \tr (x)=0\}$$ 
by putting 
$$P(\sum_{i=1}^{3}x_ie_i)=P(x_1,x_2,x_3)$$
for an orthonormal basis $\{ e_i\}$ of $D_\infty^{(0)}$ with
respect to the norm form $n$.
The space $U_{\nu}^{(0)} $ is known to have a basis of rational
polynomials (i.e., polynomials that take rational values on vectors in
$D^{(0)}=D_\infty^{(0)}\cap D$).
 
The group $D_\infty^\times/\R^\times$
acts on  $U_{\nu}^{(0)}$
through the
representation
$\tau_{\nu}$ 
(of highest weight $(\nu)$)
given by  $$(\tau_{\nu}(y))(P)(x)=P(y^{-1}xy).$$

Changing the orthonormal basis above  amounts
to replacing $P$ by $(\tau_{\nu}(y))(P)$ for some $y \in
D_\infty^\times.$

By $\dla\quad,\quad\dra_0$ we denote the suitably normalized
invariant scalar product in the representation space
$U_\nu^{(0)}$. 

For  $\nu_1 \ge 
\nu_2$ 
the 
$H^+(\R)$-space
$$U_{\nu_1}^{(0)}\otimes 
U_{\nu_2}^{(0)}$$ (irreducible of highest weight
$(\nu_1+\nu_2,\nu_1-\nu_2))$
is isomorphic to the $H^+(\R)$-space
$U_{\nu_1,\nu_2}$ of $\C[X_1, X_2]$-valued harmonic  
forms on $D_{\infty}^2$ transforming according to the
representation of 
$GL_2(\R)$ of highest weight $(\nu_1+\nu_2,\nu_1-\nu_2).$

An intertwining map $\Psi$ has been given explicitly in \cite[Section 
3]{bs-nachrichten}; for ${\bf x}=(x_1,x_2)\in D_{\infty}^2$ the
polynomial $\Psi(Q)({\bf x})\in \C[X_1, X_2]$ is homogeneous of degree
$2\nu_2.$  We write now for $Q \in  U_{\nu_1}^{(0)}\otimes 
U_{\nu_2}^{(0)}$
\begin{equation}
  \label{eq:psizerlegung}
  \Psi(Q)({\bf x})=\Sum_{\alpha_1+\alpha_2=2\nu_2}c_{\alpha_1 \alpha_2}({\bf
  x},Q)X_1^{\alpha_1} X_2^{\alpha_2}.
\end{equation}
The map ${\bf x} \mapsto c_{\alpha_1 \alpha_2}({\bf
  x},Q)$ is (for fixed $Q$) a polynomial in $x_1,x_2$ that is harmonic
  of degree $\alpha_1'=\alpha_1+\nu_1-\nu_2$ in $x_1$ and harmonic of degree
  $\alpha_2'=\alpha_2+\nu_1-\nu_2$ in $x_2,$ 
  and for $h\in H^+({\bf R})$ we have 
$$c_{\alpha_1 \alpha_2}(h{\bf
  x},Q)= c_{\alpha_1 \alpha_2}({\bf
  x},h^{-1}Q).$$ 

The irreducibility of the space $ U_{\nu_1}^{(0)}\otimes 
U_{\nu_2}^{(0)}$ implies that this map is nonzero for some $Q.$
We 
denote by $U_{\alpha}$ the space of harmonic
  polynomials of degree $\alpha$ on  $D_{\infty}$ with invariant 
scalar product $\dla \quad,\quad \dra$. If $\alpha$ is even, the
$H(\R)$-spaces $U_\alpha$ and $U_{\alpha/2}^{(0)}\otimes
U_{\alpha/2}^{(0)}$ are isomorphic and will be identified.
 
The map 
\begin{equation}
  \label{eq:trildef1}
(Q, R_1, R_2)\mapsto \dla c_{\alpha_1 \alpha_2}(
  \cdot,Q),R_1\otimes R_2 \dra
\end{equation}
for $Q \in  U_{\nu_1}^{(0)}\otimes 
U_{\nu_2}^{(0)}, R_1 \in U_{\alpha_1'}, R_2\in U_{\alpha_2'}$
defines
  then a nontrivial invariant trilinear form for the triple of
  $H({\bf   R})$-spaces $((U_{\nu_1}^{(0)}\otimes 
U_{\nu_2}^{(0)}), U_{\alpha_1'}, U_{\alpha_2'}).$ 

\begin{lemma}\label{existtrilin}
Let integers $\nu_1\ge \nu_2$ and 
$\beta_1,\beta_2$  be given for which  
$$\beta_1'=\beta_1+\nu_1-\nu_2,\beta_2'=\beta_2+\nu_1-\nu_2$$ are even.  
Then there exists a nontrivial
$H({\bf R})$- invariant trilinear form $T$ on the space 
$\Unn \otimes U_{\beta_1'} \otimes U_{\beta_2'}$ if and only if
there exist integers $\alpha_1,\alpha_2,\gamma$ such that
$\beta_i=\alpha_i+\gamma$ and $\alpha_1+\alpha_2=2\nu_2$
holds. This form is unique up to scalar multiples and can be
decomposed as
$$T=T_1^{(0)}\otimes 
 T_2^{(0)}$$ with (up to scalars) unique nontrivial
 invariant trilinear forms 
$$T_i^{(0)}=T_{i,\beta_1',\beta_2'}^{(0)}$$  on 
$$U^{(0)}_{\nu_i}\otimes U^{(0)}_{\beta_1'/2}\otimes
U^{(0)}_{\beta_2'/2}. $$ 
In particular, for $\gamma=0$ and $T$ fixed, the
trilinear form given in (\ref{eq:trildef1}) is proportional to $T$
(with a nonzero factor
$\tilde{c}(\nu_1,\nu_2,\alpha_1,\alpha_2))$.   
\end{lemma}
{\it Proof.}
Decomposing  
$$  U_{\beta_1'}=U^{(0)}_{\beta_1'/2}\otimes U^{(0)}_{\beta_1'/2},
 U_{\beta_2'}=U^{(0)}_{\beta_2'/2}\otimes
 U^{(0)}_{\beta_2'/2}$$ as a 
$D_\infty^\times/\R^\times\times
D_\infty^\times/\R^\times$-space one sees that $T$ as
asserted exists if and only if there are
nontrivial
invariant trilinear forms
$$T_i^{(0)}=T_{i,\beta_1',\beta_2'}^{(0)}$$  on 
$$U^{(0)}_{\nu_i}\otimes U^{(0)}_{\beta_1'/2}\otimes
U^{(0)}_{\beta_2'/2}$$ for $i=1,2.$ 
In this case $T$ decomposes as 
$$T=T_1^{(0)}\otimes 
 T_2^{(0)}.$$
The $ T_i^{(0)}$ are known to exist if and only if
the triples $$(\beta'_1/2,\beta'_2/2,\nu_1),
 (\beta'_1/2,\beta'_2/2,\nu_2)$$ are balanced, i.e, the numbers in
 either triple are the lengths of the sides of a triangle
(and then the form is unique up to scalars); they are unique
up to scalar multiplication. It is then easily checked that  
the numerical condition given above is equivalent to 
the existence of nonnegative integers $\alpha_1,\alpha_2,\gamma$
satisfying
$\beta_i=\alpha_i+\gamma$ and $\alpha_1+\alpha_2=2\nu_2.$

Consider now the Gegenbauer polynomial $G^{(\alpha)}(x,x')=$ obtained from 
$$G^{(\alpha)}_{1}
(t) \ = \ 2^{\alpha }
\ \mathop\sum\limits^{[{\alpha\over 2}]}_{j=0} (-1)^j \ 
{{1}\over{j!(\alpha- 2j)!}} \ 
{{(\alpha-j)!}\over 2^{2j}} \ t^{\alpha-2j}$$ 
by
$$ \tilde G^{(\alpha)}(x,x')= 2^\alpha
(n(x)n(x'))^{\alpha/2}G_1^{(\alpha)}({{{\tr} (x  
\overline{x'})}\over {2\sqrt{n(x)n(x')}}})$$
and normalize the scalar product on $U_{\alpha}$  such that
$G^{(\alpha)}$ is a reproducing kernel, i.\ e.\,  
$$\dla G^{(\alpha)}(x,x'),Q(x)\dra_{\alpha}
=Q(x')$$ 
for all $Q \in  U_{\alpha}$. 
Then for  $\alpha_1,\alpha_2,\alpha_1',\alpha_2'$ as above and some
fixed $Q\in U_{\nu_1,\nu_2}$ the map 
$$(x_1,x_2)\mapsto
T(Q,G^{(\alpha'_1)}(x_1,\cdot),G^{(\alpha_2')}(x_2,\cdot))$$
defines  a
polynomial $R_Q(x_1,x_2)$ in $x_1,x_2$ that is harmonic   of degree
$\alpha_1'=\alpha_1+\nu_1-\nu_2$ in $x_1$ and harmonic of degree
$\alpha_2'=\alpha_2+\nu_1-\nu_2$ in $x_2,$  
  and for $h\in H^+({\bf R})$ we have $R_Q(h{\bf
  x})= R_{h^{-1}Q}({\bf
  x}).$ 

As above we can therefore conclude that one has
\begin{equation}
  \label{eq:trilinequality}
 c_{\alpha_1 \alpha_2}({\bf
  x},Q) =\tilde{c}(\nu_1,\nu_2,\alpha_1,\alpha_2)
  T(Q,G^{(\alpha_1')}(x_1,\cdot),G^{(\alpha_2')}(x_2,\cdot),
\end{equation}
where the factor of proportionality
$\tilde{c}(\nu_1,\nu_2,\alpha_1,\alpha_2)$ is not zero.

We denote by
$$\cA(D^\times_{\A},R^\times_{\A},\nu)$$
the space of
functions $\varphi: D^\times_{\A}
\to U_\nu^{(0)}$ satisfying 
$\varphi(\gamma x u)=\tau_\nu(u^{-1}_\infty)\varphi(x)$
for $\gamma\in D^\times_{\Q}$ and $u=u_\infty u_f\in
R^\times_{\A}$, 
where $$R^\times_{\A} = D^\times_\infty\times \prod_p
R^\times_p$$ is the adelic group of units of $R$. 
These functions are determined by their values on the
representatives $y_i$ of a double coset decomposition
$$D^\times_{\A}=\cup^r_{i=1}D^\times y_i R^\times_{\A}$$
(where we choose the $y_i$ to satisfy  $y_{i,\infty} = 1$ and
$n(y_i)=1).$ 

The natural inner product on the space
$\cA(D^\times_{\A},R^\times_{\A},\nu)$  
is given by 
$$\langle\varphi,\psi\rangle=\sum_{i=1}^r\frac{\dla\varphi(y_i),
\psi(y_i)\dra_0}{e_i},$$
where $e_i=\mid(y_i Ry_i^{-1})^\times\mid$ is the number of units of
the order $R_i = y_i Ry_i^{-1}$ of $D.$

On the 
space $\cA(D^\times_{\A},R^\times_{\A},\nu)$ we have for 
$p{\hskip 1pt\not\vert} N$ 
(hermitian) Hecke operators $\tilde T(p)$ (given explicitly by the
$\Endo(U_\nu^{(0)})$-valued 
Brandt matrices $(B_{ij}(p))$)  
and for $p \mid N$ involutions $\widetilde{ w_p}$ commuting with the
Hecke operators and with each other.
 
For $i=1,2$ and $\nu_1\ge \nu_2$ with $\nu_1-\nu_2$ even we consider
now functions $\varphi_i$
in $\cA(D^\times_{\A},R^\times_{\A},\nu_i).$

\medskip
The Yoshida lifting (of degree $2$) of the pair
 $(\varphi_1,\varphi_2)$ is then given as
    \begin{multline}
      \label{eq:y2}
Y^{(2)}(\varphi_1,\varphi_2)(Z)(X_1,X_2)=\\
=\mathop\sum\limits^r_{i,j{=}1} {1\over{e_i e_j}} 
\mathop\sum\limits_{{{(x_1,x_2)}}\in (y_i Ry^{-1}_j)^2} \Psi( 
\varphi_1(y_i)\otimes \varphi_2(y_j))(x_1,x_2)(X_1,X_2)\times\\
\times \exp(2\pi i{\rm tr }
\left(\begin{pmatrix} n(x_1)& tr(\overline{x_1}x_2)\\
 tr(\overline{x_1}x_2)&n(x_2)
\end{pmatrix}Z\right).
\end{multline}

\medskip
This is a vector valued holomorphic Siegel modular form for the group
$\Gamma_0^{(2)}(N)$ with trivial character and with respect
to the representation $\sigma_{2\nu_2}\otimes \det^{\nu_1-\nu_2+2}$
(where   $\sigma_{2\nu_2}$ denotes the $2\nu_2$-th symmetric power
representation of $GL_2$). 

\medskip
If we consider the restriction  of such a modular
form to the diagonal 
$\left(\begin{smallmatrix}
z_1 &0\\0&z_2
\end{smallmatrix}\right),$ the coefficient of $X_1^{\alpha_1}
X_2^{\alpha_2}$ becomes a function 
$F^{(\alpha_1,\alpha_2)}(z_1,z_2)$
which is in both variables a scalar valued modular form
for the group
$\Gamma_0(N)$ with trivial 
character of weight $$\alpha_1+\nu_1-\nu_2+2 \text{ in } z_1,
\alpha_2+\nu_1-\nu_2+2 \text{ in } z_2.$$ 
In particular the
weights in the variables $z_1,z_2$ add up to $2\nu_1+4$
for each pair $(\alpha_1,\alpha_2)$ with $\alpha_1+\alpha_2=2\nu_2$ and the
coefficient of  $X_1^{\alpha_1} X_2^{\alpha_2}$ vanishes unless
$$\alpha_1'=\alpha_1+\nu_1-\nu_2,
\alpha_2'=\alpha_2+\nu_1-\nu_2$$ are even.

\medskip
If $f_1,f_2$ are elliptic modular forms of weights $k_1,k_2$
we define then 
$$\langle 
F\left(\left(
\begin{smallmatrix}z_1&0\\0&z_2
\end{smallmatrix}
\right)\right),f_1(z_1) f_2(z_2)\rangle_{k_1,k_2}$$
to be the double Petersson product
$$\langle \langle
F^{(\alpha_1,\alpha_2)}(z_1,z_2),f_1(z_1)\rangle_{k_1},
f_2(z_2)\rangle_{k_2}$$ if 
$$k_1=\alpha_1+\nu_1-\nu_2+2, k_2=\alpha_2+\nu_1-\nu_2+2$$
for some $\alpha_1,\alpha_2$ with $\alpha_1+\alpha_2=2\nu_2$
and to be zero otherwise (this definition coincides with the
Petersson product of the corresponding  automorphic
forms on the groups $Sp_2(\A)$ or $Sp_2(\R)$ (restricted to the naturally
embedded
$SL_2(\A)\times SL_2(\A)$) and $SL_2(\A)$ (or the respective real groups)).

\medskip
We will mainly consider Yoshida liftings for pairs of forms that 
are eigenforms of all Hecke operators and of all the involutions.
It is then easy to see that $Y^{(2)}(\varphi_1,\varphi_2)(Z)$ is
identically zero unless $\varphi_1,\varphi_2$ have the same eigenvalue
under the involution $\widetilde{w_p}$ for all $p\mid N$. The precise
conditions under which the lifting is  nonzero have been stated
in \cite{bs-nagoya2}.

\smallskip
We will finally need some facts about the correspondence 
studied e.g. in \cite{E,H-S,Sh,J-L} between
modular forms for $\Gamma_0(N)$ (with trivial character) and
automorphic forms on the adelic quaternion algebra $D^\times_{\A}.$

We consider the essential part $$\cA_{\text{ess}}
(D^\times_{\A},R^\times_{\A},\nu)$$ 
consisting of functions $\varphi$
that are orthogonal to all
$\psi\in\cA(D^\times_{\A},(R'_{\A})^\times,\nu)$ 
for orders $R'$ strictly containing $R$; this space  is invariant under
the $\tilde T(p)$ 
for $p{\hskip 1pt\not\vert} N$ and the 
$\tilde {w}_p$ for $p\mid N$
and hence has a basis of common eigenfunctions of all the
$\tilde T(p)$ for 
$p{\hskip 1pt\not\vert}N$ and all the involutions $\tilde w_p$ for
$p\mid N$.  Being the components of eigenvectors of a rational matrix
with real eigenvalues the values of these  eigenfunctions are real,
(i.e., polynomials with real coefficients in the vector values case)
when suitably normalized. 

Moreover the eigenfunctions 
are in one to one correspondence with the
newforms in the space $$S^{2+2\nu}(N)$$ of elliptic cusp forms of weight
$2+2\nu$ for the group $\Gamma_0(N)$ that are eigenfunctions
of all Hecke 
operators (if $\tau$ is the trivial representation and $R$
is a maximal order one has to restrict here to functions
orthogonal to the constant function $1$ on the quaternion
side in order to obtain cusp forms on the modular forms
side). This correspondence (Eichler's correspondence)  
preserves Hecke eigenvalues for 
$p{\hskip 1pt\not\vert}N$, and if $\varphi$ corresponds to
$f\in S^{2+2\nu}(N)$ then the eigenvalue of $f$ under the Atkin-Lehner
involution $w_p$ is equal to that of $\varphi$ under $\tilde w_p$
if $D$ splits at $p$ and 
equal to minus that of $\varphi$ under $\tilde w_p$ if
$D_p$ is a skew field. 
The correspondence can be explicitly described by associating to
$\varphi$ the modular form
$$h(z)=\mathop\sum\limits^r_{i,j{=}1} {1\over{e_i e_j}} 
\mathop\sum\limits_{{{x}}\in (y_i Ry^{-1}_j)}  
(\varphi(y_i)\otimes \varphi(y_j))(x)\exp(2\pi i n(x)z)$$
(where as above $\varphi(y_i)\otimes \varphi(y_j)$ denotes the
harmonic polynomial in $U_{2\nu}$ obtained by identifying
$U_\nu^{(0)}\otimes U_\nu^{(0)}$ with $U_{2\nu}.$)

An extension of Eichler's
correspondence  to forms $\vp$ as above that are not
essential but eigenfunctions of all the involutions
$\widetilde{w_p}$ has been given in \cite{hashi,bs-triple}.

\section{Computation of periods}

Our goal is the computation of the periods 
$$\langle 
Y^{(2)}(\varphi_1,\varphi_2)\left(\left(
\begin{smallmatrix}z_1&0\\0&z_2
\end{smallmatrix}
\right)\right),f_1(z_1) f_2(z_2)\rangle_{k_1,k_2}$$
defined above 
for elliptic modular forms $f_1,f_2$ for the group
$\Gamma_0(N).
$
For this  we study first how the vanishing of this period integral
depends on the eigenvalues of the functions involved under the
Atkin-Lehner involutions or their quaternionic and Siegel modular
forms counterparts.

For a Siegel modular form $F$ for the group $\Gamma_0^{(2)}(N)$ we
let the Atkin-Lehner involutions with respect to the variables
$z_1,z_2$ act on the restriction of $F$ to the diagonal matrices $\left(
\begin{smallmatrix}z_1&0\\0&z_2
\end{smallmatrix}
\right)$ and 
denote by $F\left(\left(
\begin{smallmatrix}z_1&0\\0&z_2
\end{smallmatrix}
\right)\right)|\widetilde{W_p}$ the result of this action.

\begin{lemma}

Let $N$ be squarefree and let $F$ be a vector valued Siegel
modular form of degree $2$ 
for the representation $\rho$ of $GL_2({\C})$ of highest
weight $(\lambda_1,\lambda_2)$ in the space
$\C[X_1,X_2]_{\lambda_1-\lambda_2}$ of homogeneous
polynomials of degree ${\lambda_1-\lambda_2}$ in $X_1,X_2$
with respect to $\Gamma^{(2)}_0(N)$ and assume for $p\mid N$
that the restriction of $F$ to the diagonal is an
eigenform of $\widetilde{W_p}$ with eigenvalue $\tilde{\epsilon}^{(0)}_p$.  
Let $f_1,f_2$
be elliptic cusp forms of weights $k_1,k_2$ for
$\Gamma_0(N)$ that are eigenforms of the Atkin-Lehner
involution $w_p$ with eigenvalues
$\epsilon_p^{(1)},\epsilon_p^{(2)}$. 
Then the period integral 
\begin{equation}
\langle 
F\left(\left(
\begin{smallmatrix}z_1&0\\0&z_2
\end{smallmatrix}
\right)\right),f_1(z_1) f_2(z_2)\rangle_{k_1,k_2}
\end{equation}
is zero unless one has 
$\tilde{\epsilon_p^{(0)}}\epsilon_p^{(1)}\epsilon_p^{(2)}=1.$

\end{lemma}
\vskip0.3cm
{\it Proof.}
Applying the Atkin-Lehner involution $w_p$ to both variables
$z_1,z_2$ one sees that this is obvious.

\smallskip
We can view the condition of Lemma 2.1 as a (necessary)
local condition for the nonvanishing of the period integral
at the finite primes dividing the level, with a similar
role being played at the infinite primes by the condition
that modular forms of the weights $k_1,k_2$ of $f_1,f_2$ appear in the
decomposition of the 
restriction of the vector valued modular form $F$ to the
diagonal (or of a suitable form in the representation space
of $F$, see below).

\begin{lemma}\label{zero_period_lemma}
Let $f_1,f_2,h_1,h_2$ be modular forms for $\Gamma_0(N)$
that are  eigenfunctions
of all Atkin-Lehner involutions for the $p\mid N$ with
$f_1,f_2$ cuspidal. Let $h_1,h_2$
have the same
eigenvalue $\epsilon'_p$ for all the $w_p$ for the $p\mid N $ and
denote by $\epsilon_p^{(1)},\epsilon_p^{(2)}$ the Atkin-Lehner
eigenvalues at $p\mid N$ of $f_1,f_2.$

For a factorization $N=N_1N_2$ where
$N_1$ has an odd number of prime factors let $D_{N_1}$ be
the quaternion algebra over $\Q$ that is ramified precisely
at $\infty$ and the primes $p\mid N_1$ and $R_{N_1}$ an
Eichler order of level $N$ in $D_{N_1}.$ 

Let
$\varphi_1^{(N_1)},\varphi_2^{(N_1)}$ be the forms in
$\cA((D_{N_1})_{\A}^\times, (R_{N_1})_{\A}^\times, \tau_i)
(i=1,2)$ corresponding to $h_1,h_2$ under Eichler's
correspondence.

Then the period integral
\begin{equation}
\langle 
Y^{(2)}(\varphi_1^{(N_1)},\varphi_2^{(N_1)})\left(
\begin{pmatrix}z_1&0\\0&z_2
\end{pmatrix}
\right),f_1(z_1) f_2(z_2)\rangle_{k_1,k_2}
\end{equation}
is zero unless $\epsilon'_p\epsilon_p^{(1)}\epsilon_p^{(2)}=-1$ holds
for precisely those $p$ that divide $N_1$; in particular it is always
zero unless $\prod_{p\mid
  N}\epsilon'_p\epsilon_p^{(1)}\epsilon_p^{(2)}=-1$ holds. 
\end{lemma}
{\it Proof.}
For each factorization of $N$ as above we denote by
$\widetilde{\epsilon_p}(N_1)$ the eigenvalue under
$\widetilde{w_p}$ of $\varphi_1^{(N_1)},\varphi_2^{(N_1)},$
we have $\widetilde{\epsilon_p}(N_1)=-\epsilon_p'$ for the
$p$ dividing $N_1$ and
$\widetilde{\epsilon_p}(N_1)=\epsilon_p'$ for the $p$
dividing $N_2$. Hence the product
$\widetilde{\epsilon_p}(N_1)\epsilon_p^{(1)}\epsilon_p^{(2)}$
is $1$ for all $p$ dividing $N$ if 
$N_1$ is the product of the primes $p\mid N$ such that
$\epsilon'_p\epsilon_p^{(1)}\epsilon_p^{(2)}=-1$ and is $-1$
for at least one $p\mid N$ otherwise; in particular a decomposition
for which
$\widetilde{\epsilon_p}(N_1)\epsilon_p^{(1)}\epsilon_p^{(2)}=1$ for
all $p\mid N$ holds and $N_1$ has an odd number of prime factors
exists if and only if we have $\prod_{p\mid
  N}\widetilde{\epsilon_p}(N_1)\epsilon_p^{(1)}\epsilon_p^{(2)}=-1.$ 

The $\widetilde{W_p}$-eigenvalue of the restriction of
$Y^{(2)}(\varphi_1^{(N_1)},\varphi_2^{(N_1)})$ to the
diagonal is  $\widetilde{\epsilon_p}(N_1)$ by the result of Lemma 9.1 of
\cite{bs-nagoya1} on the eigenvalue of
$Y^{(2)}(\varphi_1^{(N_1)},\varphi_2^{(N_1)})$ under the analogue for
Siegel modular forms of the Atkin-Lehner involution. The assertion
then follows from the previous lemma.

\smallskip

For simplicity we will in the sequel assume that $h_1,h_2,f_1,f_2$ are
all newforms 
of (square free) level $N$; essentially the same results can
be obtained  for more general quadruples of  forms of square free
level using the methods of \cite{bs-triple}.

\begin{lemma}
  Let $N\ne 1$ be squarefree, $D,R$ as described in Section 1, let
  $f_1,f_2$ be normalized newforms of weights 
  $k_1,k_2$ for the group $\Gamma_0(N)$. Let $\varphi_1,\varphi_2 \in
  \cA(D_{\A}^\times, R_{\A}^\times,\tau_i)$ be as above.
Assume that the (even) weights $k_1,k_2$ of $f_1,f_2$ can be
written as  $k_i=\alpha_i+\nu_1-\nu_2+2=\alpha_i'+2$ with
nonnegative integers $\alpha_i$   satisfying $\alpha_1+\alpha_2=2\nu_2
$ and denote for $i=1,2$ by  $\psi_i$   the $U_{\alpha_i'/2}$-valued
form in $\cA(D_{\A}^\times,   R_{\A}^\times,\tau_{\alpha'_i})$
corresponding to $f_i$ under Eichler's correspondence.
 
Then the period integral
\begin{equation}
  \label{eq:period2}
\langle Y^{(2)}(\varphi_1,\varphi_2)\left(
\begin{pmatrix}z_1&0\\0&z_2
\end{pmatrix}
\right),f_1(z_1) f_2(z_2)\rangle_{k_1,k_2}  
\end{equation}
has the (real) value 
\begin{multline}
  \label{eq:period2b}
  c\langle f_1,f_1\rangle \langle f_2,f_2   \rangle(\sum_{j=1}^r
     T_{\nu_1,\alpha_1',\alpha_2'}^{(0)}(\varphi_1(y_i)\otimes
     \psi_1(y_i)\otimes\psi_2(y_i)))\\ 
\times(\sum_{j=1}^rT_{\nu_2,\alpha_1',\alpha_2'}^{(0)} (\varphi_2(y_i)\otimes \psi_1(y_i)\otimes\psi_2(y_i))),
\end{multline}
with a nonzero constant $c$ depending only on $\nu_1,\nu_2,k_1,k_2.$
\end{lemma}

{\em Proof.}
The coefficient of $X_1^{\alpha_1}X_2^{\alpha_2}$ of the $(i,j)$-term
in $F(z_1,z_2)=Y^{(2)}(\varphi_1,\varphi_2)\left(
\begin{smallmatrix}z_1&0\\0&z_2
\end{smallmatrix}
\right)$ is (with $Q_{ij}:=\varphi_1(y_i)\otimes\varphi_2(y_j)$) equal to 

\begin{multline}\label{trilinyosh}
\tilde{c}(\nu_1,\nu_2,\alpha_1,\alpha_2)
  \sum_{(x_1,x_2)\in I_{ij}}T(Q_{ij},G^{(\alpha_1')}(x_1,\cdot),G^{(\alpha_2')}(x_2,\cdot))\\
  \times\exp(2\pi in(x_1)z_1)\exp(2\pi i n(x_2)z_2)
\end{multline}
 by
(\ref{eq:trilinequality}). We write
\begin{equation}
\Theta_{ij}^{(\alpha')} (z)(x')=\sum_{x\in
  I_{ij}}G^{(\alpha')}(x,x')\exp(2\pi in(x)z)
\end{equation}
for the $U^{\alpha'}$-valued theta series attached to $I_{ij}$ and the
Gegenbauer polynomial $G^{(\alpha)}(x,x')$
and rewrite (\ref{trilinyosh}) as
\begin{equation}
  \label{trilinyosh2}
  \tilde{c}(\nu_1,\nu_2,\alpha_1,\alpha_2)
  T(Q_{ij},\Theta_{ij}^{(\alpha_1')}(z_1),\Theta_{ij}^{(\alpha_2')}(z_2)).
\end{equation}
The $ij$-term of the period integral (\ref{eq:period2})
becomes then 
\begin{equation}
  \label{period3}
  \tilde{c}(\nu_1,\nu_2,\alpha_1,\alpha_2)
  T(Q_{ij},\langle\Theta_{ij}^{(\alpha_1')}(z_1), f_1(z_1)\rangle,
  \langle \Theta_{ij}^{(\alpha_2')}(z_2), f_2(z_2)\rangle,
\end{equation}
which by
  (3.13) of \cite{bs-hamburg} and the factorization 
  \begin{equation}
    \label{eq:trilinfactor}
    T=T_{\nu_1,\nu_2,\alpha_1',\alpha_2'}=T_{\nu_1,\alpha_1',\alpha_2'}^{(0)}
    \otimes T_{\nu_2,\alpha_1',\alpha_2'}^{(0)}=T_1^{(0)}\otimes T_2^{(0)}
  \end{equation}
is equal to

\begin{multline}  \label{eq:ij-partresult}  
\tilde{c}(\nu_1,\nu_2,\alpha_1,\alpha_2)\langle f_1,f_1\rangle
  \langle f_2,f_2 \rangle
  T_1^{(0)}(\varphi_1(y_i),\psi_1(y_i),\psi_2(y_i))
\\\times T_2^{(0)}(\varphi_2(y_j),\psi_1(y_j),\psi_2(y_j)).
\end{multline}

Summation over $i,j$ proves the assertion. The value computed is real
since the values of the $\vp_i,\psi_i$ are so and since $T_0$ is known
to be real.

\begin{theorem}\label{maintheorem}
  Let $h_1,h_2,f_1,f_2,\psi_1,\psi_2$ be as in Lemma 2.2,
  Lemma 2.3 with Atkin-Lehner eigenvalues $\epsilon_p'$ for $h_1,h_2$
  and $\epsilon_p^{(1)},\epsilon_p^{(2)}$ for $f_1,f_2$; assume 
  $\prod_{p\mid 
  N}\epsilon'_p\epsilon_p^{(1)}\epsilon_p^{(2)}=-1$.
Let $D$ be  the quaternion algebra over $\Q$ which is ramified
  precisely at the primes $p\mid N$ for which
  $\epsilon'_p\epsilon_p^{(1)}\epsilon_p^{(2)}=-1$ holds and $R$ an
  Eichler order of level $N$ in $D,$ let $\varphi_1,\varphi_2$ be the
  forms in 
 $\cA(D_{\A}^\times, R_{\A}^\times,\tau_{1,2})$ corresponding to
  $h_1,h_2$ under Eichler's correspondence.

Then the square of the period integral 
\begin{equation}
  \label{eq:period4}
\langle Y^{(2)}(\varphi_1,\varphi_2)\left(
\begin{pmatrix}z_1&0\\0&z_2
\end{pmatrix}
\right),f_1(z_1) f_2(z_2)\rangle_{k_1,k_2}    
\end{equation}
is equal to
\begin{equation}
  \label{eq:result}
  \frac{c}{\langle h_1,h_1 \rangle \langle h_2,h_2 \rangle} L(h_1,f_1,f_2;\frac{1}{2}) L(h_2,f_1,f_2;\frac{1}{2}),
\end{equation}
where $c$ is an explicitly computable  nonzero number
depending only on $\nu_1,\nu_2,k_1,k_2, 
N$ and the triple product $L$-function $L(h,g,f;s)$ is normalized to
have its functional equation under $s\mapsto 1-s.$

In particular the period integral is nonzero if and only if the
central critical value of $L(h_1,f_1,f_2;s)L(h_2,f_1,f_2;s)$ is nonzero.
\end{theorem}

{\em Proof.} The choice of the decomposition $N=N_1N_2$ made above
implies that we can use Theorem 5.7 of \cite{bs-triple} to express
the right hand side of (\ref{eq:period2b}) by the product of central
critical values of the triple product $L$-functions associated to
$(h_1,f_1,f_2),(h_2,f_1,f_2).$ The Petersson norms of $f_1,f_2$ appearing in
Theorem 5.7 of \cite{bs-triple} cancel against those appearing in the
proof of Lemma 2.3.

\medskip
\begin{remark} \begin{itemize}
\item[a)] If the product  $\prod_{p\mid 
  N}\epsilon'_p\epsilon_p^{(1)}\epsilon_p^{(2)}$ is $+1$ we know from
\cite{bs-triple} that the sign in the functional equation of the triple
product $L$-functions  $L(h_1,f_1,f_2;s),L(h_2,f_1,f_2;s)$ is $-1$ and
hence the central critical values are zero; from Lemma 2.2 we know
that for any Yoshida
lifting $F$ associated to $h_1,h_2$ as in Lemma 2.2
the Petersson product of the restriction of $F$ to the diagonal and
$f_1(z_1),f_2(z_2)$ is zero as well. 
\item[b)] It should be noticed that given $h_1,h_2$ there
  are $2^{\omega(N)-1}$ possible choices of the quaternion
  algebra with respect to which one considers the Yoshida
  lifting associated to $h_1,h_2.$ All these Yoshida
  liftings are different, but have the same Satake
  parameters for all $p\nmid N.$ Given $f_1,f_2$ with
$\prod_{p\mid 
  N}\epsilon'_p\epsilon_p^{(1)}\epsilon_p^{(2)}=-1$ there is
then precisely one choice of quaternion algebra that leads
to a nontrivial result for the period integral, all the
others give automatically zero by Lemma \ref{zero_period_lemma}. The choice 
of this quaternion algebra should be
seen as variation of the Vogan $L$-packet of the $p$-adic
component of the adelic representation generated by the Siegel
modular form for the $p \mid N$ in such a way that the
resulting $L$-packet satisfies the local Gross-Prasad
condition for the split $5$-dimensional and $4$-dimensional
orthogonal groups.
\end{itemize}
\end{remark}
\medskip
We want to rephrase the result of Theorem \ref{maintheorem} in order to replace the
factor of comparison $\langle h_1,h_1 \rangle \langle
h_2,h_2 \rangle$ occurring by a factor depending only on 
$F=Y^{(2)}(\varphi_1,\varphi_2)$ instead of $h_1,h_2.$
Concerning the symmetric square $L$-function of $F$
occurring in the following corollary we remind the reader
that we view $F$ as as an automorphic form on the adelic
orthogonal group of the $5$-dimensional quadratic space $V$ of
discriminant $1$ over ${\bf Q}$ that contains a 2-dimensional
totally isotropic subspace.

\begin{corollary}\label{cor_syml2}
Under the assumptions of  Theorem \ref{maintheorem} and the additional
assumption that $h_1,h_2$ are not proportional, the value of
(\ref{eq:result}) is equal to:  
\begin{equation}
  \label{corsyml}
\frac{c \langle F,F \rangle}{L^{(N)}(F,\text{Sym}^2,1)}
L(h_1,f_1,f_2;\frac{1}{2}) L(h_2,f_1,f_2;\frac{1}{2}),     
\end{equation}
where again $c$ is an explicitly computable nonzero constant depending
only on the 
levels and weights involved and $L^{(N)}(F,\text{Sym}^2,s)$
is the $N$-free part of the $L$ function of $F$ with respect
to the symmetric square of the $4$ dimensional
representation of the $L$-group of the group $SO(V)$ ($V$ as
above).  
\end{corollary}

{\em Proof.} Since $h_1,h_2$ are not proportional, the Siegel modular form $F$
is cuspidal and the Petersson product  $\langle F,F \rangle$ is well
defined. From \cite[Proposition 10.2]{bs-nagoya1} we
recall that $\langle F,F \rangle$ is (up to a nonzero
constant) equal to the the 
residue at $s=1$ of the $N$-free part $D_F^{(N)}(s)$ of the
degree $5$ $L$-function  associated to $F$ (normalizing the
$\varphi_i$ to $\langle \varphi_i,\varphi_i \rangle=1$; the
formulas given in \cite{bs-nagoya1} generalize easily to the
situation where the $\varphi_i$ take values in
harmonic polynomials). It is also well known that $\langle
h_i,h_i \rangle$ is equal (up to a nonzero constant
depending only on weights and levels) to
$D_{h_i}^{(N)}(1)$ where $D_{h_i}^{(N)}(s)$ is the symmetric
square $L$-function associated to $h_i$. Comparing the
parameters of the $L$-functions  $L^{(N)}(F,\text{Sym}^2,s)$
and $D_F^{(N)}(s)D_{h_1}^{(N)}(s)D_{h_2}^{(N)}(s)$ we see
that the value of  $L^{(N)}(F,\text{Sym}^2,s)$ at $s=1$ is
equal to the residue at $s=1$ of
$D_F^{(N)}(s)D_{h_1}^{(N)}(s)D_{h_2}^{(N)}(s),$ which gives
the assertion.

\medskip

{\em Remark. }
We can as well view $L^{(N)}(F,\text{Sym}^2,s)$
as the exterior square of the degree $5$ $L$-function
associated to $F$.

\medskip
Let us discuss now two degenerate cases:
\begin{corollary}\label{degenerate_cases}
  \begin{itemize}
  \item[a)] Under the assumptions of Theorem \ref{maintheorem} replace
    $\varphi_2$ by the constant function $(\sum_i
      \frac{1}{e_i})^{-1}$  (and hence $h_2(z)$ by the
    Eisenstein series $E(z)=(\sum_{i,j}
    \frac{1}{e_ie_j})^{-1}\sum_{i,j}
    \frac{1}{e_ie_j}\Theta_{ij}^{(0)}(z)).$ Then  
    the period integral 
\begin{equation}
  \label{eq:period4b}
\langle Y^{(2)}(\varphi_1,\varphi_2)\left(
\begin{pmatrix}z_1&0\\0&z_2
\end{pmatrix}
\right),f_1(z_1) f_2(z_2)\rangle_{k_1,k_2}    
\end{equation}
is zero unless $f_1=f_2=:f,$ in which case its square
is equal to the value at $s=1$ of 
\begin{equation}
  \label{degenerate_syml2}
\frac{c \langle F,F \rangle}{L^{(N)}(F,\text{Sym}^2,s)}
L(h_1,f,f;s-\frac{1}{2}) L(E,f,f;s-\frac{1}{2}),     
\end{equation}
where again $c$ is an explicitly computable nonzero constant depending
only on the 
levels and weights involved and $L(E,f_1,f_2;s)$ is defined in the same way as
the triple product $L$-function for a triple of cusp forms, setting the $p$-parameters
of $E$ equal to $p^{1/2},p^{-1/2}$ for $p\nmid N.$
\item[b)]   Under the assumptions of Theorem \ref{maintheorem} let $h=h_1=h_2.$
Then the period integral (\ref{eq:period4}) is equal to
\begin{equation}
  \frac{c}{\langle h,h \rangle } L(h,f_1,f_2;\frac{1}{2}),
\end{equation}
where $c$ is a nonzero constant depending only on the
levels and weights involved.

\end{itemize}

\end{corollary}
{\em Proof. }
Both assertions are obtained in the same way as Theorem 2.4 and
Corollary \ref{cor_syml2}; 
notice that in case a) (with $\omega(N)$ denoting the number of prime factors of
$N)$ both  
$L(E,f,f;s-\frac{1}{2})$ and $L^{(N)}(F,\text{Sym}^2,s) $ are of order
$\omega(N)-1$ at $s=1.$ 

{\em Remark}
\begin{itemize}
\item[a)] The form of our result given in Corollary \ref{cor_syml2}
  could in principle be true for any Siegel modular form $F$
  instead of a Yoshida lifting if one replaces
  $L(h_1,f_1,f_2;\frac{1}{2}) L(h_2,f_1,f_2;\frac{1}{2})$ by
  the value  $L({\mathrm Spin}(F),f_1,f_2,\frac{1}{2})$ of the
  spin $L$-function mentioned
  in the introduction. There is, however, not much known
  about the analytic properties of
  $L^{(N)}(F,\text{Sym}^2,s),$ in particular this $L$-
  function might have a zero or a pole
  at $s=1.$
\item[b)] In the degenerate case of Corollary \ref{degenerate_cases}
  a) the Yoshida lifting $F$ is 
 the Saito-Kurokawa lifting associated to $h_1.$ The result of that case could
 also be true in the case that $h_1,f_1,f_2$ are  of level $1$ and $F$ is the
 Saito-Kurokawa lifting of $h_1,$ but we can not prove this at present (except
 for the vanishing of the period integral in the case that $f_1 \ne f_2$, which
 is easily proved).

 We notice that in that last case (as well as in the related case of a Yoshida
 lifting of Saito-Kurokawa type) the period integral is seen to be equal to
 the Petersson product $\langle \phi_1(\tau, 0), f(\tau) \rangle,$ where 
$\phi_1(\tau, z)$ is the first Fourier Jacobi coefficient of $F.$

In the case of Corollary \ref{degenerate_cases}  b) the Yoshida lift
$F$ can be viewed as an 
Eisenstein series of Klingen type associated to $h,$ in particular its image
under Siegel's $\Phi$-operator is equal to $h$ (more precisely, the Klingen
Eisenstein series in question is a sum of Yoshida liftings associated to various
quaternion algebras of level dividing $N$ (see \cite{bs-hamburg}), where the
other contributions yield a vanishing period integral). One checks that the result of
Corollary \ref{degenerate_cases}  b) is valid also for $F$ denoting
the Klingen Eisenstein  series 
attached to a cuspidal normalized Hecke eigenform $h$ of level $1$ and $f_1,f_2$
two cuspidal normalized Hecke eigenforms of level $1.$
 
\end{itemize}

\medskip
We can obtain a result similar to that of Theorem
\ref{maintheorem} for more general weights $k_1,k_2$ of the
modular forms $f_1,f_2.$ For this, remember that according
to Lemma \ref{existtrilin} the value given in
(\ref{eq:period2b}) for the period integral in question also
makes sense if one replaces throughout $\alpha'_1,\alpha'_2$
by $\beta'_1=\alpha'_1+\gamma, \beta'_2=\alpha'_2+\gamma$
for some fixed $\gamma>0;$
the forms $f_1,f_2$ then having weights $k_i=\alpha'_i+2+\gamma$
for $i=1,2.$ As noticed above, our period integral becomes
$0$ in this situation. We can, however, modify the function  
$Y^{(2)}(\varphi_1,\varphi_2)(Z)$ by a differential
operator $\tdiff$ in such a way that 
$\tdiff Y^{(2)}(\varphi_1,\varphi_2)$ is a function on
${\bf H}\times {\bf H}$ 
that is a modular form of 
weights $k_1,k_2$ of $z_1,z_2$ as described above and yields a
value for the 
period integral of the same form as the one given in form
(\ref{eq:period2b}).  

More precisely, we have:

\begin {proposition} \label{diffprop} For nonnegative integers $k$, $r$ and $l$ with 
$k\geq 2$ and any partition $l=a+b$,
there exists a (non-zero) holomorphic differential operator
${\mathcal D}^r_{k,a,b}$ (polynomial in $X_1^2{\partial\over \partial z_{1}},
X_1X_2{\partial\over \partial z_{12}}, 
X_2^2{\partial\over \partial z_{2}}$, evaluated in $z_{12}=0$) mapping
${\bf C}[X_1,X_2]_l$-valued functions on ${\bf H_2}$ to 
${\bf C}\cdot X_1^{a+r}X_2^{b+r}$-valued functions
on ${\bf H}\times {\bf H}$ and satisfying
$$ {\mathcal D}_{k,a,b}^{r}\left(F\mid_{k,l}
M_1^{\uparrow}
M_2^{\downarrow}\right)
=({\mathcal D}_{k,a,b}^rF)\mid_{k+a+r}^{z_1}
M_1\mid_{k+b+r}^{z_2}
M_2
$$ 
for all  $M_1, \, M_2 \in SL(2,{\bf R})$; here the upper indices $z_1$ 
and $z_2$ at the slash-operator indicate the variable, with respect to which
one has to apply the elements of $SL(2,{\bf R})$ and $\uparrow\downarrow $
denote the  standard embedding of $SL(2)\times SL(2)$ into $Sp(2)$ given by
$$\left(\begin{array}{cc} a & b \\ c & d \end{array}\right)^{\uparrow}\times
\left(\begin{array}{cc} A & B \\ C & D\end{array}\right)^{\downarrow}
=
\left(\begin{array}{cccc}
a & 0 & b & 0\\
0 & A & 0 & B \\
c & 0 & d & 0 \\
0 & C & 0 & D \end{array}\right)$$
Of course one can consider ${\mathcal D}_{k,a,b}^rF$
as a ${\bf C}$-valued function.

\end{proposition}
\
\begin {remark} \label{uniqueoperator} One may indeed show, that there
  exists (up to multiplication  
by a constant) precisely one such nontrivial holomorphic differential
operator.\ 
\end{remark}
\begin {corollary} The differential operator 
${\mathcal D}^r_{k,a,b}$ defined above
gives rise to a map 
$$ \tilde{{\mathcal D}}^r_{k,a,b}:M_{k,l}^2(\Gamma^2_0(N))\longrightarrow 
M^1_{k+a+r}(\Gamma^1_0(N))\otimes M^1_{k+b+r}(\Gamma^1_0(N))$$
of spaces of modular forms. (It is easy to see that for $r>0$ this map
actually goes into spaces of cusp foms).\end{corollary}
\begin {corollary} \label{pluriharmonic_polynomial} 
Denoting by $Sym_2({\bf C})$ the space of complex symmetric matrices
of size 2
we define a polynomial function
$$Q:Sym_2({\bf C})\longrightarrow {\bf C}[X_1,X_2]_{2r}$$
by
$${\mathcal D}^r_{k,a,b}e^{tr(TZ)}=Q(T)e^{t_1z_1+t_2z_2}$$
where $Z=\left(\begin{array}{cc} z_1 & z_3\\ z_3 &z_2\end{array}\right)\in
{\bf H}_2$.\\
Furthermore we assume (with $k={m\over 2}+\nu$) that 
$P:{\bf C}^{(m,2)}\longrightarrow {\bf C}[X_1,X_2]_l$
is a polynomial function satisfying\\
\begin{itemize}
\item[a)] $P$ is pluriharmonic
\item[b)] $P((X_1,X_2)A)=\rho_{\nu,l}(A)P(X_1,X_2)$ for all 
$A\in GL(2,{\bf C})$.
\end{itemize}
Then, for ${\bf Y}_1, {\bf Y_2}\in {\bf C}^m$
$$({\bf Y}_1,{\bf Y}_2)\longmapsto \left\{
P({\bf Y}_1,{\bf Y}_2)\cdot Q(\left(
\begin{array}{cc}
{\bf Y}_1^t{\bf Y}_1 & {\bf Y}_1^t{\bf Y}_2 \\
{\bf Y}_2^t{\bf Y}_1 & {\bf Y}_2^t{\bf Y}_2 \end{array}\right) )
\right\}_{a+r+\nu,b+r+\nu}$$
defines an element of $H_{a+r+\nu}(m)\otimes H_{b+r+\nu}(m)$, where
$H_{\mu}(m)$ is the space of harmonic polynomials in $m$ variables (for  the
standard quadratic form),
homogeneous of degree $\mu$ and for any $R\in {\bf C}[X_1,X_2]_{l+2r}$
we denote by $\{R\}_{\alpha,\beta}$ the coefficient of 
$X_1^{\alpha}X_2^{\beta}$ in $R$,
$\alpha+\beta=l+2r$.\end{corollary} 
The proof of Corollary \ref{pluriharmonic_polynomial} is a
vector-valued variant of similar statements in 
\cite{fourierII} and \cite[p.200]{bs-squares}, using the proposition 
above and the characterization of harmonic polynomials by 
the Gau{\ss}-transform; we leave the details of proof to the
reader.

\medskip
{\it Proof of Proposition \ref{diffprop}.}\\
We start from a Maa{\ss}-type differential operator $\delta_{k+l}$
which maps ${\bf C}[X_1,X_2]_l$-valued functions on ${\bf H}_2$ to
${\bf C}[X_1,X_2]_{l+2}$-valued ones and satisfies
$$\left(\delta_{k+l}F\right)\mid_{k,l+2}M=
\delta_{k+l}\left(F\mid_{k,l}M\right)$$
for all $M\in Sp(2,{\bf R})$.\\
It is well known, how such operators arise from elements of the universal
enveloping algebra of the complexified Lie algebra of $Sp(2,{\bf R})$, 
see e.g. \cite{harris} .
In our case (we refer to \cite{bsy} for details) we can 
describe these operators quite explicitly 
in terms of the simple operators
$$DF:= \left({1\over 2\pi i} {\partial\over \partial Z}\right)[{\bf X}]$$
and
$$NF:= \left(-{1\over 4\pi} (ImZ)^{-1}F\right)[{\bf X}]$$
Here ${\bf X}$ stands for the column vector 
$\left(\begin{array}{c} X_1\\ X_2 \end{array}\right)$.
Then we define
$$\delta_{k}F= kNF+DF$$
It is remarkable (and already incorporated in our notation!) that 
$\delta_{k+l}$ depends only on $k+l$.\\
The iteration
$$\delta_{k+l}^r:= 
\delta_{k+l+2r-2}\circ\cdots\circ\delta_{k+l+2}\circ\delta_{k+l}$$
can also be described explicitly by
$$\delta_{k+l}^r=\sum_{i=0}^r {\Gamma(k+l+r)\over \Gamma(k+l+r-i)}
{r \choose i } N^i D^{r-i}$$

For a function $F:{\bf H}_2\longrightarrow {\bf C}[X_1,X_2]_l$
and a decomposition $l=a+b$ we put
$$\nabla^r_{k+l}(a,b)F=:X_1^{a+r}X_2^{b+r}-\mbox{coefficient of} \,\left(
\delta_{2,k+l}^r F\right)_{\mid {\bf H}\times {\bf H}}$$
Then $\nabla$ has already the transfomation properties 
required in the proposition, i.e.
$$\nabla^r_{k+l}(a,b) \left(F\mid_{k,l}
M_1^{\uparrow}
M_2^{\downarrow}\right)
= \left(\nabla^r_{k+l}(a,b)F\right)
\mid_{k+a+r}^{z_1}
M_1\mid_{k+b+r}^{z_2}
M_2
$$ 
for $M_1,\,M_2\in SL(2,{\bf R})$.\\  
Moreover, if $F$ is in addition a holomorphic function on ${\bf H}_2$, then
$\nabla^r_{k+l}(a,b)F$ is a nearly holomorphic function in the sense
of Shimura (with respect to both variables $z_1$ and $z_2$), as polynomials in 
${1\over y_1}$ and ${1\over y_2}$ they are of degree $\leq r  $.
Shimura's structure theorem on nearly holomorphic functions \cite{Sh76}
says that all nearly holomorphic functions on ${\bf H}$ can be 
obtained from holomorphic functions by applying Maa{\ss} type operators
$$\delta_k:={k\over 2 iy}+{\partial\over \partial z}$$
and their iterates. This however is only true if the weight (i.e.
$k+a+r$ or $k+r+b$) is bigger than $2r$, which is not necessarily true
in our situation. We therefore use a weaker version
of Shimura' theorem (see \cite[Theorem 3.3]{Sh94}), valid under the
assumption "$w>1+r$", where $w$ is the weight at hand and $r$ is the 
degree of the nearly holomorphic function: 
Every 
such function $f$ on ${\bf H}$ of degree $\leq r$
has an expression
$$f=f_{hol}+L_w(\tilde{f})$$
where $f_{hol}$ is holomorphic and $\tilde{f}$ is again nearly holomorphic
of degree $\leq r$; in this expression 
$$L_w:= \delta_{w-2}\left(y^2{\partial\over \partial\bar{z}}\right)=
{w\over 2i}y{\partial\over \partial\bar{z}}+ 
y^2{\partial^2\over \partial z\partial\bar{z}}$$
is a "Laplacian" of weight $w$
commuting with the $\mid_w$-action of $SL(2,{\bf R})$.
We also point out that $f_{hol}$ is uniquely  determined by $f$ (in particular,
$f=f_{hol}$, if $f$ is holomorphic) and we have
$\left(f\mid_wM\right)_{hol}=\left(f_{hol}\right)\mid_wM$ 
for all $M\in SL(2,{\bf R})$.\\
If we apply this statement to $\nabla^r_{k+l}(a,b)F$, considered as function of 
$z_1$ and $z_2$, we get an expression of type
$$ \nabla^r_{k+l}(a,b)F=f+L_{k+a+r}^{z_1}g_1+
L_{k+r+b}^{z_2}g_2+L_{k+a+r}^{z_1}L_{k+r+b}^{z_2}h$$
where $f,g_1,g_2,h$ are nearly holomorphic functions on 
${\bf H}\times {\bf H}$,
$f$ being holomorphic in both variables, $g_1$ holomorphic in $z_2$,
$g_2$ holomorphic in $z_1$. Note that (due to our assumption $k\geq 2$ )
Shimura's theorem is applicable here. An inspection of Shimura's proof
(which is quite elementary for our case) shows that $f$ is indeed of the form
$f={\mathbf D} F$, where ${\mathbf D}$ is a
polynomial $p$ in
${\partial\over \partial z_{1}},
{\partial\over \partial z_{12}}, 
{\partial\over \partial z_{2}}$, evaluated in $z_{12}=0$.
This polynomial does not depend on $F$ at all and it has the 
required transformation properties. \\
It remains however to show that 
${\mathbf D}$ is not zero:\\
For this purpose, we consider the special function
$$ {\bf z}_{12}^r:
\left\{
\begin{array}{ccc}
{\bf H}_2 & \longrightarrow  &{ \bf C}[X_1,X_2]_l \\
Z & \longmapsto & z_{12}^rX_1^aX_2^b 
\end{array}\right.$$
It is easy to see that
$\nabla^r_{k+l}(a,b)({\bf z}_{12}^r)$ is then equal to the constant 
function $r!$, therefore
$$\nabla^r_{k+l}(a,b)({\bf z}_{12}^r)={\mathbf D}({\bf z}_{12}^r)=r!,$$
in particular, ${\mathbf D}$ is non-zero and
we may put
$${\mathcal D}^r_{k,a,b}=p(X_1^2{\partial\over \partial z_{1}},
X_1X_2{\partial\over \partial z_{12}}, 
X_2^2{\partial\over \partial z_{2}}),$$
evaluated at $z_{12}=0$

\medskip
We can then prove in the same way as above:

\begin{corollary}
The assertions of Lemma 2.3 and Theorem 2.4 remain true if 
$f_1,f_2$ have weights $k_i=\alpha'_i+2+\gamma (i=1,2)$ with some
$\gamma >0,$ if one replaces $$Y^{(2)}(\varphi_1,\varphi_2)\left(
\begin{pmatrix}z_1&0\\0&z_2
\end{pmatrix}
\right) $$
by
$$\tdiff Y^{(2)}(\varphi_1,\varphi_2)\left(
z_1,z_2 \right)$$
\end{corollary}

\begin{remark}\label{rep_remark} 
Application of the differential operator to \linebreak
$Y^{(2)}(\varphi_1,\varphi_2)\left(
Z)
\right)$ before restriction to the diagonal does not change
the $Sp_2(\R)$-representation space of that function, i.e.,
we have found a different function in the same
representation space whose period integral assumes the value
that is predicted by the conjecture of Gross and
Prasad. 
More precisely, (using remark \ref{uniqueoperator} and some additional
considerations) one can 
show that the vanishing of this predicted value is
already sufficient for the vanishing of the period integral
for all triples 
$F', f'_1,f'_2$  
of functions in the Harish-Chandra
modules generated by
the original functions $F,f_1,f_2$.
To obtain a similar statement for the local
representations at the finite places not dividing the level
one would have to show that a nonvanishing  invariant linear
functional 
on the tensor product of the representations is not zero on
the product of the spherical (or class 1) vectors invariant
under the maximal compact subgroup. This is expected to be true
as well; we plan to come back to these problems in future work. 

\end{remark}

\section{Restriction to an embedded Hilbert modular surface}

To avoid technical difficulties we deal here only with the simplest
case: The quaternion algebra $D$ is ramified at all primes
$p$ dividing the level $N$ and we have $\nu_1=\nu_2=0$, i.e., the
Yoshida lifting is a scalar valued Siegel modular form of
weight $2$ and the order $R$ we are considering is a maximal
order. We put $F={\bf Q}(\sqrt{N})$ and assume that $N$ is
such that the class number of $F$ is 1. We denote by $\Delta$ the
discriminant of $F,$ by $a\mapsto a^\sigma$ its nontrivial
automorphism and consider the
basis ${1,w}$ with $w=\frac{\Delta+\sqrt{\Delta}}{2}$ of the ring
${\frak o}_F $ of $F$. Denoting by $C$
the matrix 
\begin{equation*}
C:=\bigl(\begin{matrix}
1 &1\\
w&\bar{w}
\end{matrix}
\bigr)
\end{equation*}
we have the usual modular embedding
\begin{equation}
  \label{eq:modularembed}
   \iota:(z_1,z_2) \mapsto  C 
  \begin{pmatrix}
    z_1 &0\\
    0 &z_2
  \end{pmatrix}{}^t C 
\end{equation}of  ${\mathbf H}\times {\mathbf H}$ into the Siegel upper half
    plane ${\mathbf H}_2$ and 
\begin{equation} 
    \tilde{\iota}:\begin{pmatrix}
    a&b\\
    c&d
  \end{pmatrix}
\mapsto 
\begin{pmatrix}
  C&0\\
  0&{}^tC^{-1}
\end{pmatrix}
\begin{pmatrix}
  a&0&b&0\\
0&a^\sigma&0&b^\sigma\\
c&0&d&0\\
0&c^\sigma&0&d^\sigma
\end{pmatrix}\begin{pmatrix}
C^{-1}&0\\
0 &{}^tC
\end{pmatrix}
\end{equation}
from $SL_2(F)$ into $Sp_2({\mathbf Q}).$

We have then for $\gamma\in SL_2(F):$ 
\begin{equation}
  \label{eq:commute}
\tilde{\iota}(\gamma)(\iota((z_1,z_2))=\tilde{\iota}(\gamma((z_1,z_2)))  
\end{equation}
with the usual actions of the groups $SL_2(F)$ on ${\mathbf H}\times
{\mathbf H}$ and of $Sp_2({\mathbf Q})$ on the Siegel upper half
plane.

We put now 
\begin{equation*}
J=\begin{pmatrix}  
1&0&0&0\\
0&0&0&1\\
0&0&1&0\\
0&1&0&0
\end{pmatrix}
\end{equation*}
and consider for a Siegel modular form $f$ of weight $k$ for
the group $\Gamma \subseteq Sp_2(\bQ)$ the
function 
\begin{equation}
  \label{eq:pullback}
  \tilde{f}(\tau_1,\tau_2):=f|_k J (\iota(\tau_1,\tau_2)).
\end{equation}

Writing
$\iota_0=J \circ \iota, \widetilde{\iota_0}(\gamma):=J
\iota(\gamma) J^{-1}$
we see that $\tilde{f}$ is a Hilbert modular form for the
group $\widetilde{\iota_0}^{-1}(\Gamma).$

By calculating $\widetilde{\iota_0}(\gamma)$ explicitly for
$\gamma\in SL_2(F)$ one checks that $
\widetilde{\iota_0}(\gamma)$ is in $\Gamma_0^{(2)}(N)$ if
and only if $\gamma$ is in $SL_2({\mathfrak o}_F\oplus {\mathfrak
  d}) $, the group of matrices $\left(\begin{smallmatrix} a&b\\ c&d
\end{smallmatrix}\right)$ with $a,d \in {\mathfrak o}_F, c \in
{\mathfrak d}, b \in {\mathfrak d}^{-1}$, where ${\mathfrak
  d}$ is the different of $F.$

If $L$ is a ${\mathbf Z}$-lattice of (even) rank $m=2k$ with
quadratic form $q$ and associated bilinear form
$B(x,y):=q(x+y)-q(x)-q(y)$ satisfying $q(L) \subseteq
 {\bZ}$
and $N q(L^\#)\Z ={\bZ}$ ($N$ is the level of $(L,q)$) then it
is shown in \cite{bs-nagoya1} that
\begin{equation}
  \label{eq:thetaseries}
  \vartheta^{(2)}(L,q,Z)|_k J=c_{1}\sum_{x_1 \in L, x_2\in
  L^\#} \exp\bigl(2\pi i {\rm tr} \bigl(\begin{pmatrix}
  q(x_1)&B(x_1,x_2)/2\\
B(x_1,x_2)/2&q(x_2)
\end{pmatrix} Z\bigr)\bigr),
\end{equation}
where $c_{1}$ is a nonzero constant depending only on the
genus of $(L,q)$ and where $\vartheta^{(2)}(L,q,Z)$ is the
usual theta series of 
degree $2$ of $(L,q).$
One checks therefore that, writing $K=\{x_1+x_2 w \in L
\otimes F \mid x_1 \in L, x_2 \in L^\# \},$ we have 
\begin{equation}
  \label{eq:thetarestrict}
  \vartheta^{(2)}(L,q,\widetilde{\iota_0}(z_1,z_2))=\vartheta (K,q,(z_1,z_2)),
\end{equation}
where we denote by $\vartheta (K,q,(z_1,z_2))=\sum_{y \in
  K}\exp(2\pi i (q(y)z_1+q(y)^\sigma z_2))$ the theta series
  of the ${\mathfrak o}_F$-lattice $K$ with the extended
  form $q$ on it.

It is again easily checked that $L^\# \subseteq N^{-1}L$
implies that $K$ is an integral unimodular ${\mathfrak
  o}_F$-lattice, and it is well known that then the theta
series $\vartheta (K,q,(z_1,z_2))$ is a modular form of
weight $k$ for the group $SL_2({\mathfrak o}_F\oplus {\mathfrak
  d}).$

\begin{lemma}\label{jlcorresp}
Let $\tilde{D}$ be a quaternion algebra over $F$ ramified at
both infinite primes and let $\tilde{R}$ be a maximal order
in $\tilde{D}.$ Let
$$\cA(\tilde{D}^\times_{\A},\tilde{R}^\times_{\A},0)=:\cA(\tilde{D}^\times_{\A},\tilde{R}^\times_{\A})$$
be defined in the same way as in Section 1 for $D$ and let 
$\cA(\tilde{D}^\times_{\A},\tilde{R}^\times_{\A})$ be
equipped with the natural action of Hecke operators
$T({\mathfrak p})$ for the ${\mathfrak p}$ not dividing $N$
described by Brandt matrices as explained  in
\cite{eichler}.Then by associating to a Hecke eigenform
$\psi \in \cA(\tilde{D}^\times_{\A},\tilde{R}^\times_{\A})$
the Hilbert modular form 
$$f(z_1,z_2)=\int_{(\tilde{D}^\times\backslash
  \tilde{D}_{\A}^\times)\times (\tilde{D}^\times\backslash
  \tilde{D}_{\A}^\times)}  \psi(y)\psi(y')
\vartheta(y'\tilde{R}y^{-1},(z_1,z_2))dy \,dy'$$ one gets a
bijective correspondence between the $\psi$ as above and the
Hecke eigenforms of weight $2$ and trivial character for the
group $\Gamma_0({\mathfrak n},{\mathfrak d})$ of precise
level ${\mathfrak n}$ giving an explicit realization of the
correspondence of Shimizu und Jacquet/Langlands
\cite{Sh,J-L}. Here ${\mathfrak n}$ denotes the product of
the prime ideals ramified in $\tilde{D}$ and
$\Gamma_0({\mathfrak n},{\mathfrak d})$ is the subgroup of
$SL_2({\mathfrak o}_F\oplus {\mathfrak d})$ whose lower left
entries are in ${\mathfrak n}{\mathfrak d}.$ 

The function $\rho(y,y')$ on $\tilde{D}^\times_\A \times
\tilde{D}^\times_\A$ given by setting $\rho(y,y')$ equal to
the Petersson product of $f$ with
$\vartheta(y'\tilde{R}y^{-1},(z_1,z_2))$ is proportional to
$(y,y')\mapsto \psi(y)\psi(y').$

\end{lemma}
\vskip0.3cm
{\it Proof.}
The first part of this Lemma is due to Shimizu
\cite{Sh} (taking into account that by \cite{eichler} the
group $\Gamma_0({\mathfrak n},{\mathfrak d})$ is the correct
transformation group for the theta series in question). Let
$\hat{\rho}$ be the function on the adelic 
orthogonal group of $\tilde{D}$ induced by $\rho$ and let
$\hat{\psi} $ be the  function on the adelic
orthogonal group of $\tilde{D}$ induced by $(y,y')\mapsto
\psi(y)\psi(y').$ The function $\hat{\psi}$ generates an
irreducible representation space of $\tilde{D}^\times_{\A}$
whose theta lifting to $SL_2(F_{\A})$ is generated by $f,$
and $\hat{\rho}$ is a vector in the theta lifting of this
latter representation of $SL_2(F_{\A}),$ which by
\cite{moeglin-base} coincides with the original
representation space 
generated by $\hat{\psi}.$ Since both
$\hat{\rho},\hat{\psi}$ are invariant under the same maximal
compact subgroup of $\tilde{D}_\A^\times,$ the uniqueness of
such a vector implies that they must coincide up to
proportionality. That $\hat{\rho}$ is not zero follows from
the obvious fact that $f$ by its construction can not be
orthogonal to all the theta series. 

\begin{lemma}\label{hilbertperiod}
With the above notations let  $\varphi_1, \varphi_2$
in $\cA(D^\times_{\A},R^\times_{\A},0)$ be Hecke eigenforms
with the same eigenvalue under the involutions
$\widetilde{w_p}$ for the $p\mid N$ with associated newforms
$h_1,h_2$ of weight $2$ and level $N.$ Let $f$ be a Hilbert
modular form of weight $2$ for the group $SL_2({\mathfrak
  o}_F\oplus {\mathfrak d})$ that corresponds in the way described
in Lemma \ref{jlcorresp} to the function $\psi \in
\cA(\tilde{D}^\times_{\A},\tilde{R}^\times_{\A})$ for
$\tilde{D}=D\otimes F$ and $\tilde{R}$ being the maximal
order in $\tilde{R}$ containing $R.$
Then the value of the period integral 
\begin{equation}
  \label{eq:hperiod1}
  \int_{SL_2({\mathfrak o}_F\oplus {\mathfrak d})\backslash
  H\times H}( Y^{(2)}(\varphi_1,\varphi_2)\vert_2 J (\iota((z_1,z_2)))f((z_1,z_2))dz_1dz_2
\end{equation} is equal to
\begin{equation}
  \label{eq:hperiod2}
  c_{2}\langle f,f\rangle(\sum_i \varphi_1(y_i)\psi(y_i))(\sum_i \varphi_2(y_i)\psi(y_i)), 
\end{equation}
where we identify $y_i$ with $y_i\otimes 1 \in \tilde{D}$
and where $c_{2}$ is some constant depending only on $N.$ 
\end{lemma}
In order to interpret the value obtained in
(\ref{eq:hperiod2}) in the same way as in Section 2 as the
central critical value of an $L$-function, we review briefly
the integral representation of the $L$-function that one
obtains when one replaces in a triple $(h,f_1,f_2)$ of elliptic cusp forms
the pair $(f_1,f_2)$ by one Hilbert cusp form $f$ for a real quadratic
field.

For the moment, both the Hilbert cusp form $f$ and the elliptic cusp form
$h$ can be of arbitrary even weight $k.$ 
Now we consider the Siegel type Eisenstein series of weight $k$, defined on 
${\bf H}_3$
by 
$$E^k_3(W,s)=
\sum_{\gamma =\left(\begin{array}{cc} * & * \\ C & D\end{array}\right)
\in \Gamma_0^3(N)_{\infty}\backslash\Gamma^3_0(N)}
det(CW+D)^{-k} det(\Im(\gamma<W>)^s $$
Here and in the sequel we denote by $G_{\infty}$ the subgroup of $G$
defined by "$C=0$", where $G$ is any group of symplectic matrices.\\ 
We restrict this Eisenstein series to $W=
\left(\begin{array}{cc} \tau & 0 \\ 0 & Z\end{array}\right)$ with
$\tau\in{\bf H}$, $Z\in {\bf H_2}$ and furthermore we consider then the
modular embedding with respect to $Z$.

In this way we get a function $E(\tau,z_1 ,z_2,s)$, which behaves 
like a modular
form for $\tau$ and like a Hilbert modular form for $(z_1,z_2)$ 
of weight $k$.\\
We want to compute the twofold integral
$$I(f,h,s):=
\int_{SL_2({\mathfrak o}_F\oplus \mathfrak{d})\backslash {\bf H}^2}
\int_{\Gamma_0(N)\backslash {\bf H}} 
\overline{h(\tau)f(z_1,z_2)}E(\tau,z_1,z_2,s) d\tau^*dz_1^*dz_2^* 
$$
where $dz^*=y^{k-2}dxdy$ for $z=x+iy\in {\bf H}$.\\
This can be done in several ways: One can relate this integral to similar 
ones in \cite{rallis-ps} or in \cite{garrett preprint} 
(both these works are in an adelic setting) or one can try to do it along
classical lines as in \cite{garrett annals, satoh,bs-squares}.
We sketch the latter approach here (for class number one, 
$h$ being a normalized newform of level $N$.)

The inner integration over $\tau$ (which can be done with 
$Z\in {\bf H}_2 $ instead of the embedded $(z_1,z_2)$) is the 
same as in the papers mentioned above,
producing an $L$-factor
$L_2(h,2s+2k-2)$ (with $L_2({ },{ })$  denoting the symmetric square
$L$-function) times
a Klingen type Eisenstein series 
$E_{2,1}(h,s)$, which is defined as follows:
We denote by $C_{2,1}$ the maximal parabolic subgroup of $Sp(2)$ for which
the last line is of the form $(0,0,0,*)$ and we put 
$C_{2,1}(N)=C_{2,1}({\bf Q})\cap \Gamma^2_0(N)$. Furthermore we
define a function $h_s(Z)$ on ${\bf H}_2$ by
$$h_s(Z)=h(z_1)\left({det(Y)\over y_1}\right)^s,
$$
where $z_1=x_1+iy_1$ denotes the entry in the upper left corner of $Z=X+iY\in
{\bf H}_2$. Then we put
$$E_{2,1}(h,s)(Z)=
\sum_{\gamma\in C_{2,1}\backslash\Gamma^2_0(N)}h_s(Z)\mid_k\gamma$$

 To do the second integration, one needs information 
on certain cosets: This is the only new ingredient entering the picture:\\
\begin{lemma}  A complete set of representatives for 
$ C_{2,1}(N)\backslash\Gamma^2_0(N)$
is given by
$$ \{
d(M)J^{-1}\tilde{\iota}_0(\gamma)\}
$$
with 
$\gamma$ running over 
$ SL_2({\mathfrak o}_F\oplus\mathfrak{d})_{\infty}
\backslash SL_2({\mathfrak o}_F\oplus{\mathfrak d})$,
and 
$M=\left(\begin{array}{cc}* & * \\ u & v\end{array}\right)$ 
running over those elements of 
$SL(2,{\bf Z})_{\infty}\backslash SL(2,{\bf Z})$ with $v\equiv 0(N)$, where
$d$ denotes the standard embedding of $GL(2)$ in $Sp(2)$ given by
$d(M)=\left(\begin{array}{cc} 
\left( M^{-1}\right)^{t} 
& 0 \\ 0 & M\end{array}\right)$.
\end{lemma}
This lemma is related to the double coset decomposition
$$C_{2,1}(N)\backslash \Gamma_0^2(N)/ 
\tilde{\iota}_0(SL_2({\mathfrak o}_F\oplus \mathfrak{d}))$$
and somewhat analogous to the coset decomposition in \cite[p.692]{satoh}; 
we omit the proof.\\
We may now do the usual unfolding to get
\begin{multline}
\int_{SL_2({\mathfrak o}_F\oplus \mathfrak{d})\backslash {\bf H}^2} 
\overline{h(z_1,z_2)}\widetilde{E_{2,1}(h,*,s)} (z_1,z_2)dz_1^*dz_2^*\\
\quad =\int_{SL_2({\mathfrak o}_F\oplus \mathfrak{d})_{\infty}\backslash{\bf H}^2 } 
\sum_{M=\left(\begin{array}{cc} * & * \\ v & u \end{array}\right)}
h((v,u) C
\left(\begin{array}{cc}z_1 & 0 \\
0 & z_2\end{array}\right)C^t
\left(\begin{array}{c}v \\ u \end{array}\right))\\ 
\quad\times\overline{f(z_1,z_2)}
\left( {Dy_1y_4\over (v+u\omega)^2y_1+(v+u\bar{\omega}^2y_2} \right)^s
dz_1^*dz_2^*
\end{multline}
Using the Fourier
expansions of $f$ and $F$,
$$h(z)=\sum_{n=1}^{\infty}a(n)e^{2\pi inz},
\qquad
f(z_1,z_2)=\sum_{\nu\in {\bf 0}_F,\nu\gg 0} A(\nu)e^{2\pi itr({\nu\cdot z})}$$ 
one can (after some standard calculations)
write the integral above as
$$\gamma(s)\sum_n\sum_{\sim \backslash(u,v)}a(n) 
\overline{A(n(v+u\omega)^2)}n^{-s-2k+2}N(v+u\omega)^{-2s-2k+2}$$
where we use the following equivalence relation:
two pairs $(u,v)$ and $(u',v')$ are called equivalent iff
$v+u\omega$ and $v'+u'\omega$ are equal up to a unit from ${\mathfrak o}_F$
as a factor.\\
Assume now in addition that $h$ is a normalized eigenfunction of all 
Hecke operators; then we define the $L$-function $L(h\otimes f,s)$
as an Euler product over all primes $p$ with Euler factors (at least for
$p$ coprime to $N$)
$$L_p(h\otimes f,s):=L^{Asai}_p(f,\alpha_pp^{-s})L_p^{Asai}(f,\beta_pp^{-s})$$
where we use the Euler factors $L_p^{Asai}(f,s)$ of the Asai-L-function
attached
to $f$ (see \cite{asai}) and $\alpha_p$ and $\beta_p$ are the 
Satake-p-parameters
attached to the eigenform $h$ (normalized to have absolute values
$p^{\frac{k-1}{2}}$). We will write later $L(h,f;s)$ to denote the
shift of this $L$-function that is normalized to have functional
equation under $s \mapsto 1-s$.\\  
By standard calculation, we see that the integral above is, after multiplication
by $L_2(h,2s+2k-2)$, equal to the $L$-function 
$L(h\otimes f,s+2k-2)$
(up to elementary factors; the condition $v\equiv 0(N)$ also creates some
extra contribution for $p$-Euler factors with $p\mid N$). 
This calculation of course requires some formal
calculations similar to those given e.g. in \cite{garrett annals}.\\
\begin {remark} If the class number $H$ of $F$ is different from one, 
then the orbit structure
is more complicated. One gets $H$ different sets of representatives
of the type described in the lemma above (each one twisted by a matrix in
$SL_2(F)$ mapping a cusp into $\infty$). After unfolding, one
gets then a Dirichlet series also involving Fourier coefficients of $f$
at all the $H$ different cusps.
If we assume that $f$ is the first component (i.e. the one corresponding to
the principal ideal class) in a tuple of 
$H$ Hilbert modular
forms such that the corresponding adelic modular form is an eigenform of all
Hecke operators,  
then it is possible (but quite unpleasant) to transfer that Dirichlet series
into the Euler product in question.
\end{remark}

Now we return to the case of weight $2$. We can compute the integral 
$I(f,h,s)$ at $s=0$ not only by unfolding as above but also by using the
Siegel-Weil formula for the Eisenstein series in the integrand.
Then one gets in the same way as in \cite{bs-triple}
that the square of the right hand side of
(\ref{eq:hperiod2}) is (up to an explicit constant) the
product of the central critical values  of the $L$-functions
attached to the pairs $h_1,f$ and $h_2,f$ as above:

\begin{theorem}
Let  $\varphi_1, \varphi_2,h_1,h_2,f,\psi$ be as in Lemma
\ref{hilbertperiod}. 
Then the square of the period integral 
\begin{equation}
  \int_{SL_2({\mathfrak o}_F\oplus {\mathfrak d})\backslash
  {\bf H}\times {\bf H}}
( Y^{(2)}(\varphi_1,\varphi_2)\vert_2 J (\iota((z_1,z_2)))f((z_1,z_2))dz_1dz_2
\end{equation} is equal to
\begin{equation}
\frac{c_{3}}{\langle h_1,h_1\rangle \langle h_2,h_2\rangle
  }L(h_1,f;1/2)L(h_2,f; 1/2), 
\end{equation}
where $c_{3}$ is an explicitly computable  nonzero number
depending only on $N$ and 
the product $L$-function $L(h,f;s)$ is normalized
to 
have its functional equation under $s\mapsto 1-s.$

In particular the period integral is nonzero if and only if the
central critical value of $L(h_1,f;s)L(h_2,f,s)$ is nonzero.

\end{theorem}

\vskip0.5cm
Siegfried B\"ocherer \\
Kunzenhof 4B \\
79117 Freiburg\\
Germany \\
boech@siegel.math.uni-mannheim.de\\[0.3cm]
Masaaki Furusawa\\
Department of Mathematics\\
Graduate School of Science\\
Osaka City University\\
Sugimoto 3--3--138, Sumiyoshi-ku\\
Osaka 558-8585, Japan\\
furusawa@sci.osaka-cu.ac.jp
\\[0.3cm]
Rainer Schulze-Pillot\\
Fachrichtung 6.1 Mathematik\\
Universit\"at des Saarlandes (Geb. 27)\\
Postfach 151150\\
66041 Saarbr\"ucken\\
Germany\\
schulzep@math.uni-sb.de


\newpage
   \begin{thebibliography}{99}

\bibitem{asai} T.Asai: On certain Dirichlet Series Associated with Hilbert
Modular Forms and Rankin's Method. Math.Ann.226, 81-94(1977)

\bibitem{fourierII} S. B\"ocherer: \"Uber die Fourier-Jacobi-Entwicklung
Siegelscher Eisensteinreihen II. Math. Z.189, 81-100(1985)
\bibitem{bsy} S. B\"ocherer, T. Satoh, T. Yamazaki: On the pullback 
of a holomorphic differential
operator and its application to vector-valued Eisenstein series.
Comm.Math.Univ.S.Pauli 41, 1-22(1992)

\bibitem{bs-nagoya1} S.\ B\"ocherer, R.\ Schulze-Pillot: Siegel modular
  forms and theta
series attached to quaternion algebras. Nagoya Math.J. 121(1991),
35-96
\bibitem{bs-nachrichten} S.\ B\"ocherer, R.\ Schulze-Pillot: Mellin
  transforms of  
vector valued theta series attached to quaternion algebras,
Math. Nachr. 169 (1994), 31-57  
\bibitem{bs-hamburg} S.\ B\"ocherer, R.\ Schulze-Pillot: Vector valued  theta
series and Waldspurger's theorem, Abh. Math. Sem. Hamburg 64 (1994), 211-233


\bibitem{bs-triple} S.\ B\"ocherer, R.\ Schulze-Pillot: On the central
critical value of the triple product L-function.
In: Number Theory 1993-94, 1-46. Cambridge University Press 1996


\bibitem{bs-nagoya2} S.\ B\"ocherer, R.\ Schulze-Pillot: Siegel modular
  forms and  
  theta series attached to quaternion algebras II. Nagoya
  Math.J.147(1997), 71-106 
\bibitem{bs-squares} S. B\"ocherer,R. Schulze-Pillot: Squares of automorphic
forms on quaternion algebras and central critical values of $L$-functions
of modular forms. J.Number Theory 76, 194-205(1999)

\bibitem{E} M. Eichler: The basis problem for modular forms and the traces
of the Hecke operators, p. 76-151 in Modular functions of one variable $I$,
Lecture Notes Math. 320, Berlin-Heidelberg-New York 1973
\bibitem{eichler} M. Eichler:  On theta functions of real algebraic number fields. Acta Arith. 33 (1977),
no. 3, 269--292

\bibitem{Fre 1}  E. Freitag: Siegelsche Modulfunktionen,
Berlin-Heidelberg-New York 1983

\bibitem{garrett annals} P.Garrett: Decomposition of Eisenstein series:
Rankin triple products. Annals of Math. 125, 209-235(1987) 
\bibitem{garrett preprint} P.Garrett: Integral representations of certain 
$L$-functions attached to one, two, and three modular forms.
University of Minnesota Technical report 86-131 (1986)

\bibitem{gross-kudla} B. H. Gross, S. S. Kudla: Heights and the central
  critical values of triple 
product $L$-functions. Compositio Math. 81 (1992), no. 2, 143--209
\bibitem{gropra1} B. H. Gross, D. Prasad: On the decomposition of a representation of
${\rm SO}\sb n$ when restricted to ${\rm SO}\sb {n-1}$. Canad. J. Math. 44 (1992), no. 5, 974--1002
\bibitem{gropra2} B. H. Gross, D. Prasad: On irreducible representations of ${\rm SO}\sb
{2n+1}\times{\rm SO}\sb {2m}$. Canad. J. Math. 46 (1994), no. 5, 930--950. 
\bibitem{hammond} W. Hammond: The modular groups of Hilbert
  and Siegel, Am. J.Math. 88 (1966), 497-516
\bibitem{harris} M. Harris: Special values of zeta functions attached 
to modular forms. Ann.scient.Ec.Norm.Sup.14, 77-120 (1981)  

\bibitem{HK} M. Harris, S. Kudla: The central critical value of a
triple product $L$-function, Annals of Math. 133 (1991), 605-672
\bibitem{hashi} K. Hashimoto: On Brandt matrices of Eichler orders,
 Mem. School Sci. Engrg. Waseda Univ. No.
59, (1995), 143--165 (1996)
\bibitem{H-S} H. Hijikata, H. Saito: On the representability of modular
forms by theta series, p. 13-21 in Number Theory, Algebraic Geometry and
Commutative Algebra, in honor of Y. Akizuki, Tokyo 1973

\bibitem{Ibu} T. Ibukiyama: On differential operators on automorphic forms 
and invariant pluriharmonic polynomials.Comm.Math.Univ.S.Pauli 48, 103-118(1999)
\bibitem{J-L} H. Jacquet, R. Langlands: Automorphic forms on
$GL(2)$, Lect. Notes in Math. 114, Berlin-Heidelberg-New York 1970
\bibitem{moeglin-base} C. Moeglin: Quelques propri\'etes de
  base des series 
  theta, J. of Lie Theory 7 (1997), 231-238
\bibitem{rallis-ps} I.Piatetski-Shapiro, S.Rallis: 
Rankin triple $L$-functions. Compos.Math.64, 31-115(1987)

\bibitem{prasad_seesaw} D. Prasad: Some applications of seesaw duality
  to branching laws, Math. Annalen 304, 1-20 (1996)
\bibitem{satoh} T. Satoh:  Some remarks on triple $L$-functions. Math.Ann.276,
687-698(1987) 

\bibitem{Sh} H. Shimizu: Theta series and automorphic forms on $GL_2$.
J. of the Math. Soc. of Japan 24 (1972), 638-683

\bibitem{Sh76}  G. Shimura, The special values of the zeta
functions associated with cusp forms,
Comm. pure appl. Math. 29, 783-804 (1976)

\bibitem{Sh94}  G. Shimura, Differential operators, holomorphic
projection, and singular forms, Duke Math.J.76, 141-173 (1994)

\bibitem{yoshida1} H. Yoshida: Siegel's modular forms and the arithmetic of
  quadratic forms. Invent. Math. 60 (1980), no. 3, 193--248
 \end{thebibliography}
\end{document}